\documentclass[12pt]{article}
\usepackage{amssymb}
\usepackage{amsmath}
\usepackage{amsthm}
\usepackage{enumitem}
\usepackage{graphicx}
\usepackage[total={6.5in,8.75in}, top=1.2in, left=0.9in, includefoot]{geometry}
\def\bcb{border-collision bifurcation}

\def\pwl{piecewise-linear}
\def\pws{piecewise-smooth}

\def\sw{switching manifold}

\begin{document}

\newtheorem{theorem}{Theorem}
\newtheorem{corollary}[theorem]{Corollary}
\newtheorem{lemma}[theorem]{Lemma}
\theoremstyle{definition}
\newtheorem{definition}{Definition}

\title{Simultaneous Border-Collision and Period-Doubling Bifurcations}
\author{
D.J.W.~Simpson and J.D.~Meiss\thanks{
D.~J.~W.~Simpson and J.~D.~Meiss gratefully
acknowledge support from NSF grant DMS-0707659.
}\\
Department of Applied Mathematics\\
University of Colorado\\
Boulder, CO 80309-0526}
\maketitle

\begin{abstract}
We unfold the codimension-two simultaneous occurrence
of a \bcb~and a period-doubling bifurcation
for a general \pws, continuous map.
We find that, with sufficient non-degeneracy conditions,
a locus of period-doubling bifurcations emanates non-tangentially from a locus of \bcb s.
The corresponding period-doubled solution undergoes a \bcb~along a curve
emanating from the codimension-two point and tangent to the period-doubling locus here.
In the case that the map is one-dimensional local dynamics are completely classified;
in particular, we give conditions that ensure chaos.
\end{abstract}

{\bf
Piecewise-smooth, continuous maps
are used to model nonsmooth physical situations with
discrete-time inputs, and also arise as Poincar\'{e} maps
in a variety of \pws, systems of ordinary differential equations.
Frequently a period-doubling bifurcation is the cause for
a critical change in dynamics of a map producing
a physically important period-doubled solution.
A special situation arises in a \pws, continuous map
if this bifurcation occurs at a non-differentiable point in phase space.
We show that soon after its creation,
the period-doubled solution undergoes a \bcb.
Beyond this bifurcation the period-doubled solution may persist
but the solution may also immediately become chaotic.
}

\section{Introduction}
\label{sec:INTRO}

Piecewise-smooth systems are currently being utilized in a wide variety of fields
to model physical systems involving a discontinuity or sudden change.
Examples include electrical circuits with a switching component
such as a diode or a transistor \cite{ZhMo03,BaVe01,Ts03},
vibro-impacting systems and systems with friction \cite{LeNi04,Po00,WiDe00,Br99,AwLa03},
and optimization in economics \cite{PuSu06,MoLa07,CaJa06}.
In this paper we study discrete-time, \pws~systems that are everywhere continuous.

A map
\begin{equation}
x_{n+1} = F(x_n) \;,
\label{eq:genMap}
\end{equation}
is {\em \pws~continuous} if $F : \mathbb{R^N} \to \mathbb{R^N}$ is
everywhere continuous but non-differentiable on codimension-one boundaries
called {\em \sw s}.
The collision of an invariant set with a \sw~under parameter variation of (\ref{eq:genMap})
may produce a bifurcation that cannot occur in smooth systems.
When the invariant set is a fixed point and a bifurcation occurs,
the bifurcation is called a {\em \bcb}.

Border-collision bifurcations may be analogous to familiar
smooth bifurcations, such as saddle-node and period-doubling bifurcations.
Alternatively they may be very complex.
For instance a chaotic attractor may be generated at a \bcb, even
if the map is only one-dimensional \cite{DiBu08}.
In two or more dimensions, a complete classification of
all possible dynamical behavior local to a \bcb~is yet to be determined,
see for instance \cite{ZhMo06b,SuGa06z,SiMe08b,BaGr99,DuNu99}.

At a \bcb~a fixed point, $x^*$, lies on a \sw.
Due to the lack of differentiability
the stability multipliers associated with $x^*$ are not well-defined.
However, the \sw~locally divides phase space into 
regions where the map has different smooth components;
by continuity $x^*$ is a fixed point of each component
and, restricted to each component,
associated multipliers of $x^*$ are well-defined.
Typically $x^*$ will be a hyperbolic fixed point of each component.
A special situation arises when $x^*$ is non-hyperbolic in at least
one component.

Under the assumption that the \bcb~occurs at a smooth point on a single
\sw, there are exactly two smooth map components (half-maps), in some neighborhood.
Consider the situation that $x^*$ is hyperbolic in one
half-map but non-hyperbolic in the other half-map.
Then there are three codimension-two cases to consider;
namely that on the non-hyperbolic side, $x^*$ has an associated multiplier 1,
a multiplier $-1$ or a complex conjugate pair of multipliers,
${\rm e}^{\pm 2 \pi {\rm i} \omega}$.
These correspond to the simultaneous occurrence of a \bcb~with
a saddle-node bifurcation,
a period-doubling bifurcation and
a Neimark-Sacker bifurcation, respectively.
Unlike a generic \bcb,
in each codimension-two case linear terms of the half-maps are insufficient
to describe all local dynamical behavior.

For the first case,
generically a locus of saddle-node bifurcations emanates from the codimension-two,
\bcb~which in a two-parameter bifurcation diagram is a curve
tangent to a locus of \bcb s \cite{Si08}
(see \cite{SiKo09} for the continuous-time analogue).
The third case may be very complicated and is left for future work.
The precise nature of local bifurcations may
depend on the rationality of $\omega$.
A recent analysis of these situations in a
general \pws~map is given in \cite{CoDe09}.
The purpose of the present paper is to determine generic dynamical
behavior of (\ref{eq:genMap}) local to the codimension-two \bcb~of the second case.

The coincidence of border-collision and period-doubling
has been observed in real-world systems.
Recall that a corner collision in a Filippov system provides an example
of a \bcb~in a map \cite{DiBu01c}.
In \cite{AnBe05,AnOl08} the authors describe the
coincidence of a corner collision with period-doubling
in a Filippov model of a DC/DC power converter.
Grazing-sliding (which also corresponds to a \pws, continuous Poincar\'{e} map,
see \cite{DiKo02})
and period-doubling have been seen to occur simultaneously
in a model of a forced, dry-friction oscillator \cite{KoDi06}.

In the next section we introduce a piecewise form for (\ref{eq:genMap})
suitable for the ensuing analysis and compute its fixed points.
In \S\ref{sec:UNFOLD} we unfold the codimension-two situation
for the case that the map is one-dimensional.
Here all local dynamical behavior is determined.
We find there are essentially six distinct, generic, unfolding scenarios
though each exhibit the same basic bifurcation curves.
When the map is of a dimension greater than one
we prove that the same basic bifurcation curves exist
but do not classify all possible local dynamics, \S\ref{sec:HIGHD}.
Finally \S\ref{sec:CONC} gives conclusions.

\section{Formulization of the codimension-two point and setup}
\label{sec:FORM}

We wish to explore local dynamical behavior near an arbitrary \bcb~of (\ref{eq:genMap}).
We restrict ourselves to a neighborhood 
in both phase space and parameter space of the \bcb~within which
we assume the \bcb~occurs on an isolated
$C^l$ \sw~and that away from the \sw, $F$ is $C^k$.
Then, via coordinate transformations similar to those given in \cite{DiBu01},
we may assume the \bcb~occurs at the origin, $x = 0$,
when a parameter, $\mu$, is zero,
and that to order $l$ the \sw~is simply the plane $e_1^{\sf T} x = 0$.
In this paper we are not concerned with effects due to
a nonsmooth \sw~(for studies of various \pws~systems
involving a non-differentiable \sw, see, for instance,
\cite{DiBu01c,OsDi08,BuPi06,LeVa02}).
For this reason we assume $l$ is sufficiently large
to not affect local dynamics and for simplicity assume
the \sw~is exactly the plane $e_1^{\sf T} x = 0$.

For convenience we denote the first component of the
$N$-dimensional vector, $x$, by $s$, i.e.
\begin{equation}
s = e_1^{\sf T} x \;.
\nonumber
\end{equation}
Let $\eta$ be a second map parameter, independent of $\mu$.
Then we are interested in local dynamics of
\begin{equation}
x_{n+1} = f(x_n;\mu,\eta)
= \left\{ \begin{array}{lc}
f^{(L)}(x_n;\mu,\eta), & s_n \le 0 \\
f^{(R)}(x_n;\mu,\eta), & s_n \ge 0
\end{array} \right. \;,
\label{eq:pwsMap}
\end{equation}
where
\begin{equation}
f^{(i)}(x;\mu,\eta) = \mu b(\mu,\eta) + A_i(\mu,\eta) x + O(|x|^2) + o(k) \;,
\label{eq:fiForm}
\end{equation}
is $C^k$ for $i = L,R$.
Here $b \in \mathbb{R}^N$ and $A_L$ and $A_R$ are $N \times N$ matrices
which, by the assumption of continuity of (\ref{eq:pwsMap}),
are identical in their last $N-1$ columns.
Throughout this paper we use the notation $O(k)$ [$o(k)$]
to denote terms that are order $k$ or larger [larger than order $k$]
in all variables and parameters of a given expression.
We refer to $f^{(L)}$ as the left half-map
and $f^{(R)}$ as the right half-map.

Let us first determine fixed points of (\ref{eq:pwsMap}).
Assume for the moment that for some $\eta$, 1 is not an eigenvalue of each $A_i(0,\eta)$.
Then near $(x;\mu) = (0;0)$, each half-map, $f^{(i)}$, has a unique fixed point
\begin{equation}
x^{*(i)}(\mu,\eta) = (I - A_i(0,\eta))^{-1} b(0,\eta) \mu + O(\mu^2) \;.
\label{eq:xStar}
\end{equation}
Recall that the adjugate of any matrix $X$, ${\rm adj}(X)$, obeys
${\rm adj}(X) X = \det(X) I$ \cite{Be92,Ko96}.
Since $A_L$ and $A_R$ can differ in only the first column,
it follows that the first row of
${\rm adj}(I - A_L)$ and of ${\rm adj}(I - A_R)$ are identical.
We denote this row by $\varrho^{\sf T}$:
\begin{equation}
\varrho^{\sf T} \equiv e_1^{\sf T} {\rm adj}(I - A_L) = e_1^{\sf T} {\rm adj}(I - A_R) \;.
\label{eq:varrhoDef}
\end{equation}
Multiplication of (\ref{eq:xStar}) by $e_1^{\sf T}$ on the left yields
\begin{equation}
s^{*(i)}(\mu,\eta) = \left.
\frac{\varrho^{\sf T} b}{\det(I - A_i)} \right|_{\mu = 0} \mu + O(\mu^2) \;,
\label{eq:sStar}
\end{equation}
where $s^{*(i)}$ denotes the first
component of $x^{*(i)}$.

The point $x^{*(L)}$ is a fixed point of the \pws~map
(\ref{eq:pwsMap}) and said to be {\em admissible}
(that is, a fixed point of (\ref{eq:pwsMap})) whenever $s^{*(L)} \le 0$,
otherwise it is said to be {\em virtual}.
Similarly $x^{*(R)}$ is admissible exactly when $s^{*(R)} \ge 0$.
The condition
\begin{equation}
\varrho^{\sf T} b |_{\mu = 0} \ne 0 \;,
\label{eq:rhobnonzero}
\end{equation}
ensures that the parameter $\mu$ ``unfolds'' the border-collision;
we will assume this nondegeneracy condition holds.

Let us briefly review Feigin's classification of
\bcb s by the relative admissibility of fixed points and two-cycles
(refer to \cite{DiFe99} for further details).
If $I - A_L(0,\eta)$ and $I - A_R(0,\eta)$ are both nonsingular
and (\ref{eq:rhobnonzero}) is satisfied,
for small $\mu$ each $x^{*(i)}$ is admissible
when $\mu = 0$ and for one sign of $\mu$.
Let $\sigma_i^+$ [$\sigma_i^-$] denote the number of real eigenvalues of $A_i(0,\eta)$
that are greater than 1 [less than $-1$].
By (\ref{eq:sStar}),
if $\sigma_L^+ + \sigma_R^+$ is even,
$x^{*(L)}$ and $x^{*(R)}$ are admissible for different signs of $\mu$.
In this case as $\mu$ is varied through zero
a single fixed point effectively {\em persists}.
If instead $\sigma_L^+ + \sigma_R^+$ is odd,
$x^{*(L)}$ and $x^{*(R)}$ are admissible for the same sign of $\mu$;
at $\mu = 0$ the two fixed points collide and annihilate
in a {\em nonsmooth fold bifurcation}.
Furthermore if $I - A_L(0,\eta) A_R(0,\eta)$ is nonsingular,
$-1$ is not an eigenvalue of $A_L(0,\eta)$ and $A_R(0,\eta)$,
and $\sigma_L^- + \sigma_R^-$ is odd, then a two-cycle of (\ref{eq:pwsMap})
(that collapses to the origin as $\mu \to 0$)
is admissible for one sign of $\mu$,
whereas if $\sigma_L^- + \sigma_R^-$ is even, no such two-cycle exists.

When $\mu = 0$, 
the origin is a fixed point
of (\ref{eq:pwsMap}) that lies on the \sw, $s = 0$.
For small $\mu \ne 0$, stability of, and dynamics local to, the fixed point, $x^{*(i)}$,
are determined by the eigenvalues of $A_i(0,\eta)$.
In nondegenerate situations the eigenvalues of the matrices $A_i(0,\eta)$
completely determine dynamics local to the \bcb~at $\mu = 0$.
A special situation arises whenever either $A_L(0,\eta)$ or $A_R(0,\eta)$,
has an eigenvalue that lies on the unit circle.
Without loss of generality we may assume the former matrix has an eigenvalue
on the unit circle when $\eta = 0$.
In this situation a local smooth bifurcation may occur for the left half-map,
$f^{*(L)}$, at $\mu = \eta = 0$.
The nature of the smooth bifurcation is determined by nonlinear terms,
consequently such terms are important to the \bcb.

The remainder of this paper is an analysis of the case that
$-1$ is an eigenvalue of $A_L(0,0)$.

\section{Unfolding in one dimension}
\label{sec:UNFOLD}

\begin{figure}[ht!]
\begin{center}
\includegraphics[width=15cm,height=16cm]{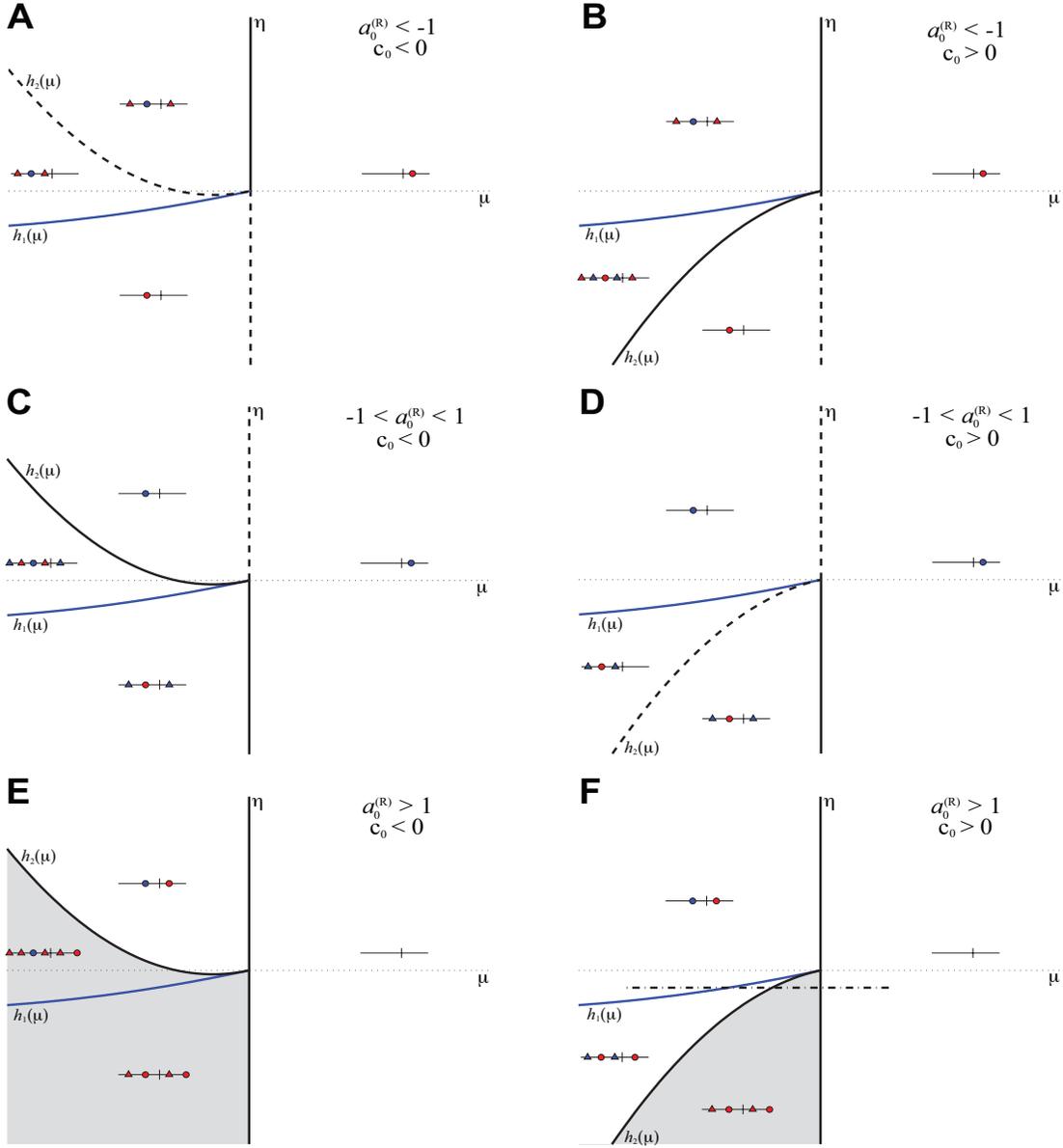}
\caption{
Sketches of bifurcation sets of (\ref{eq:1dmapPD})
for different values of $a_0^{(R)}$ and $c_0$, see Theorem \ref{th:c21dpd}.
Blue curves correspond to loci of period-doubling bifurcations;
black curves correspond to border-collisions and are shown dashed
when no topological change occurs at the border-collision.
Insets show phase portraits with circles used to denote fixed points
and triangles used to denote points on two-cycles.
These are blue when the solutions are stable and red otherwise.
The small vertical lines denote the \sw.
Chaotic dynamics occurs only for parameter values in the shaded regions;
indeed, periodic solutions of a period higher than two
only occur in these regions.
The bifurcation diagram shown in Fig.~\ref{fig:bifDiag1dex}
is taken along the dash-dot line segment of panel F.
\label{fig:bsPD}
}
\end{center}
\end{figure}

We begin by treating the one-dimensional case, that is,
$x = s \in \mathbb{R}$.
Then the map (\ref{eq:pwsMap}) may be written as
\begin{equation}
x_{n+1} = \left\{ \begin{array}{lc}
f^{(L)}(x_n;\mu,\eta), & x_n \le 0 \\
f^{(R)}(x_n;\mu,\eta), & x_n \ge 0
\end{array} \right. \;.
\label{eq:1dmapPD}
\end{equation}
We assume that the two smooth components of the map are at least $C^3$
(i.e.~$k \ge 3$) and write them as
\begin{equation}
\begin{split}
f^{(L)}(x;\mu,\eta) & =
\mu b(\mu,\eta) + a_L(\mu,\eta) x + p(\mu,\eta) x^2
+ q(\mu,\eta) x^3 + o(x^3) \;, \\
f^{(R)}(x;\mu,\eta) & =
\mu b(\mu,\eta) + a_R(\mu,\eta) x + O(x^2) \;.
\end{split}
\label{eq:1dcoeffs}
\end{equation}
The map (\ref{eq:1dmapPD}) has a \bcb~at $\mu = 0$.
We assume that the codimension-two situation of interest
occurs at $\mu = \eta = 0$, that is $a_L(0,0) = -1$.
In order for parameters $\mu$ and $\eta$ to unfold the bifurcation
therefore we require $b(0,0) \ne 0$ (by (\ref{eq:rhobnonzero}))
and $\frac{\partial a_L}{\partial \eta}(0,0) \ne 0$.
For simplicity we assume appropriate scalings upon $\mu$ and $\eta$
have been performed such that these values are both 1
as stated in the following theorem.
\begin{theorem}~\\
Consider (\ref{eq:1dmapPD}) with (\ref{eq:1dcoeffs}) and suppose
\begin{enumerate}[label=\roman{*}),ref=\roman{*}]
\item
\label{it:pdSing}
$a_L(0,0) = -1$ (singularity condition),
\item
\label{it:pdNDSW}
$b(0,0) = 1$ (border-collision non-degeneracy condition),
\item
\label{it:pdTrans}
$\frac{\partial a_L}{\partial \eta}(0,0) = 1$ (transversality condition).
\end{enumerate}
Let
\begin{eqnarray}
a_0^{(R)} & = & a_R(0,0) \;,
\label{eq:a0RPD} \\
c_0 & = & p^2(0,0) + q(0,0) \;.
\label{eq:c0PD}
\end{eqnarray}
Then there exists $\delta > 0$ such that
the left half-map restricted to the neighborhood
\begin{equation}
\mathcal{N} = \left\{ (x;\mu,\eta) ~\Big|~ |x|, |\mu|, |\eta| < \delta \right\}
\label{eq:nbd}
\end{equation}
has a unique fixed point 
given by a $C^k$ function,
$x^{*(L)}(\mu,\eta) = \frac{1}{2} \mu + O(2)$, and
there exist unique, $C^{k-1}$ functions
$h_1, h_2 : \mathbb{R} \to \mathbb{R}$ that satisfy
\begin{eqnarray}
h_1(\mu) & = & -\left( \frac{\partial a_L}{\partial \mu}
+ p \right) \bigg|_{\mu = \eta = 0} \mu + O(\mu^2) \;,
\label{eq:h1PDth} \\
h_2(\mu) & = & h_1(\mu) - \frac{c_0}{4} \mu^2 + o(\mu^2) \;,
\label{eq:h2PDth}
\end{eqnarray}
such that
\begin{enumerate}[label=\arabic{*}),ref=\arabic{*}]
\item
\label{it:resh1}
when $\mu < 0$, $x^{*(L)}(\mu,h_1(\mu))$ has an associated multiplier of $-1$,
and if $c_0 \ne 0$ the curve $\eta = h_1(\mu)$
corresponds to a locus of admissible, period-doubling bifurcations of $x^{*(L)}$,
\item
\label{it:resh2}
the origin belongs to a period-two orbit of $f^{(L)}$ on $\eta = h_2(\mu)$
that is admissible when $\mu < 0$,
\item
\label{it:resfig}
if $a_0^{(R)} \ne \pm 1$ and $c_0 \ne 0$,
then fixed points and two-cycles of (\ref{eq:1dmapPD}) exist
in the sectors shown in Fig.~\ref{fig:bsPD},
\item
\label{it:resncycles}
if $a_0^{(R)} < 1$ or $\mu > 0$ or $\eta > h_2(\mu)$,
(\ref{eq:1dmapPD}) has no $n$-cycles for any $n \ge 3$ in $\mathcal{N}$,
\item
\label{it:reschaos}
if $a_0^{(R)} > 1$ and $\mu < 0$ and $\eta < h_2(\mu)$,
then (\ref{eq:1dmapPD}) exhibits chaos in $\mathcal{N}$.
\end{enumerate}
\label{th:c21dpd}
\end{theorem}

\begin{figure}[b!]
\begin{center}
\includegraphics[width=8.4cm,height=7cm]{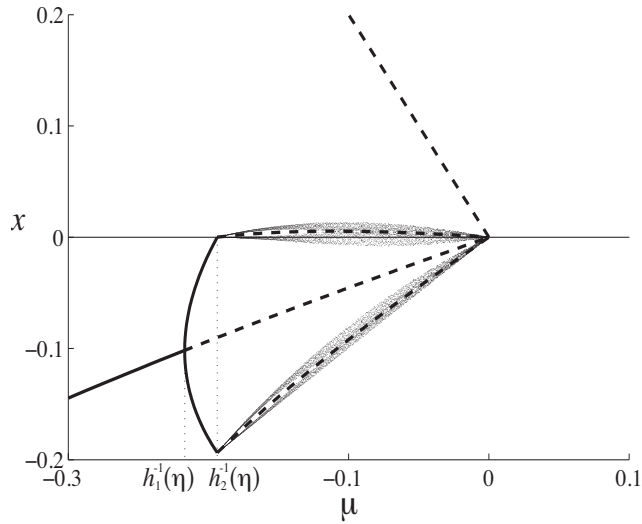}
\caption{
A bifurcation diagram for (\ref{eq:1dmapPD}) when
$b = 1$, $a_L = \eta-1$, $a_R = \frac{3}{2}$, $p = -1$, $q = \frac{3}{2}$
and there are no additional higher order terms.
The value of $\eta$ is fixed at $-0.25$.
The parameters considered correspond to the dash-dot line segment shown
in Fig.~\ref{fig:bsPD}.
Stable [unstable] fixed points and two-cycles are indicated by
solid [dashed] curves.
A supercritical period-doubling bifurcation occurs at
$\mu = h_1^{-1}(\eta) \approx -0.2169$;
a \bcb~of the period-doubled solution occurs at
$\mu = h_2^{-1}(\eta) \approx -0.1937$.
When $h_2^{-1}(\eta) < \mu < 0$ there exists a chaotic attracting set.
When $\mu > 0$ there is no local attractor.
\label{fig:bifDiag1dex}
}
\end{center}
\end{figure}

A proof of Theorem \ref{th:c21dpd} is given in Appendix \ref{sec:PROOF1D}.
To facilitate the extension of Theorem \ref{th:c21dpd} to higher dimensions
in \S\ref{sec:HIGHD},
the theorem is stated in a manner that yields a partial result when
the non-degeneracy condition on the nonlinear terms, $c_0 \ne 0$, is not satisfied.
The sign of $c_0$ determines the criticality of the period-doubling bifurcations.

Fig.~\ref{fig:bsPD} shows the six different bifurcation scenarios
depending on the value of $a_0^{(R)}$ and the sign of $c_0$
that are predicted by the theorem.
In each case the nature of the \bcb~at $\mu = 0$ for small $\eta \ne 0$
may be determined by refering to the
well-understood, one-dimensional, \pwl~map
\begin{equation}
x_{n+1} = \left\{ \begin{array}{lc}
\mu + a_L x_n, & x_n \le 0 \\
\mu + a_R x_n, & x_n \ge 0
\end{array} \right. \;,
\label{eq:1dmap}
\end{equation}
see for instance \cite{DiBu08}.
For example, along the positive $\eta$-axis in Fig.~\ref{fig:bsPD}-A,
we have $a_L \approx -1 + \eta > -1$
and $a_R < -1$, thus by Feigin's classification a single fixed point persists and
a two-cycle is created that coexists with the left fixed point
and is stable.

In panels E and F, along the negative $\eta$-axis we have
$a_L \approx -1 + \eta < -1$ and $a_R > 1$.
At points just to the left of the negative $\eta$-axis,
both fixed points and a two-cycle
are admissible and unstable.
Forward orbits are attracted to a chaotic solution
created at the \bcb~at $\mu = 0$.
Unstable high period orbits are also generated at the \bcb;
these are of even period.
To further illustrate the unfolding,
Fig.~\ref{fig:bifDiag1dex} shows a bifurcation diagram
corresponding to a one-parameter slice of Fig.~\ref{fig:bsPD}-F.
The chaotic attracting set collapses to the origin at $\mu = 0$
and to a two-cycle at $\mu = h_2^{-1}(\eta)$.

Two-cycles created in \bcb s
at $\mu = 0$ consist of one point in each half-plane.
In contrast, two-cycles created in period-doubling bifurcations
along $\eta = h_1(\mu)$ lie entirely in the left half-plane.
Along $\eta = h_2(\mu)$ these two orbits coincide
and are given explicitly by $\{ 0, \mu b(\mu,h_2(\mu)) \}$.
The second iterate of (\ref{eq:1dmapPD}) can be transformed into
a map of the form (\ref{eq:pwsMap}) along $\eta = h_2(\mu)$,
see (\ref{eq:f2proofPD}).
The resulting left and right period-two maps have
multipliers $1 - c_0 \mu^2 + o(\mu^2)$
and $-a_0^{(R)} + O(\mu)$ respectively.
By Feigin's classification,
four-cycles are born at $\eta = h_2(\mu)$
exactly when $a_0^{(R)} > 1$, i.e.~in panels E and F of Fig.~\ref{fig:bsPD},
though these are not shown.

\section{Higher dimensions}
\label{sec:HIGHD}

We now analyze the codimension-two scenario
in an arbitrary number of dimensions,
i.e.~$x \in \mathbb{R}^N$.
Recall that for smooth systems the concept of dimension reduction by
center manifold analysis allows bifurcations in systems
of any number of dimensions to be transformed to their low dimensional normal forms.
Unfortunately this technique cannot be applied here
because (\ref{eq:pwsMap}) is not
differentiable in a neighborhood of the bifurcation.
However we can compute a center manifold of the $C^k$, left half-map.
Assuming certain nondegeneracy conditions,
we will show that the restriction of the left half-map to this
center manifold yields a one-dimensional map that has the same form
as the left half-map of (\ref{eq:1dmapPD}) and satisfies the
conditions (\ref{it:pdSing})-(\ref{it:pdTrans}) of Theorem \ref{th:c21dpd}.
Consequently the aspects of Theorem \ref{th:c21dpd} that incorporate only
the left half-map extend immediately to the same scenario in higher dimensions.
That is, there exists a curve $\eta = h_1(\mu)$
along which period-doubling bifurcations occur
and a curve $\eta = h_2(\mu)$ along which the two-cycle created at the
period-doubling bifurcation collides with the \sw.
The drawback of this approach is that we fail
to describe dynamical behavior of (\ref{eq:pwsMap}) that occurs on
both sides of the \sw.

Let us state our main result:
\begin{theorem}~\\
Consider the piecewise-$C^k$, continuous map
(\ref{eq:pwsMap}) with (\ref{eq:fiForm}) on $\mathbb{R}^N$
and assume $N \ge 2$ and $k \ge 4$.
Suppose that near $(\mu,\eta) = (0,0)$,
$A_L(\mu,\eta)$ has an eigenvalue $\lambda(\mu,\eta) \in \mathbb{R}$
with an associated eigenvector, $v(\mu,\eta)$. In addition, suppose
\begin{enumerate}[label=\roman{*}),ref=\roman{*}]
\item
\label{it:pdSing2}
$\lambda(0,0) = -1$ is of algebraic multiplicity one
and $A_L(0,0)$ has no other eigenvalues on the unit circle,
\item
\label{it:pdNDSW2}
$\varrho^{\sf T}(0,0) b(0,0) \ne 0$
(where $\varrho^{\sf T}$ is given by (\ref{eq:varrhoDef})),
\item
\label{it:pdTrans2}
$\frac{\partial \lambda}{\partial \eta}(0,0) = 1$,
\item
\label{it:pdCenterMan2}
the first element of $v(0,0)$ is nonzero, thus by scaling assume
$e_1^{\sf T} v(0,0) = 1$,
\item
\label{it:pdNDRHM2}
$\det \left( I - A_L(0,0) A_R(0,0) \right) \ne 0$.
\end{enumerate}
Then, in the extended coordinate system $(x;\mu,\eta)$,
there exists a $C^{k-1}$ three-dimensional center manifold, $W^c$,
for the left half-system that passes through the origin and is not
tangent to the \sw~at this point.
Furthermore, $f^{(L)} \big|_{W^c}$, described in the coordinate system
$(\hat{x};\hat{\mu},\hat{\eta}) =
(s;\frac{2 \varrho^{\sf T} b}{\det(I-A_L)} \big|_{(0,0)} \mu,\eta)$,
is given by the left half-map of (\ref{eq:1dmapPD}) with
``hatted variables''
and the conditions (\ref{it:pdSing})-(\ref{it:pdTrans})
of Theorem \ref{th:c21dpd} will be satisfied.

If $c_0 \ne 0$ in Theorem \ref{th:c21dpd}
(where $c_0$ is computed for $f^{(L)} \big|_{W^c}$), 
then along a $C^{k-2}$ curve, $\hat{\eta} = h_1(\hat{\mu})$,
the unique fixed point of the left half-map near the origin undergoes a
period-doubling bifurcation that is admissible when $\hat{\mu} < 0$,
and along a $C^{k-2}$ curve, $\hat{\eta} = h_2(\hat{\mu})$,
the period-doubled solution collides with the \sw.
These two curves intersect and are tangent at $\hat{\mu} = \hat{\eta} = 0$.
Moreover, a unique two-cycle consisting of one point on each
side of the \sw~exists for small $\mu$ and $\eta$.
This cycle is admissible exactly when $\hat{\mu} \le 0$
and $\hat{\eta} \le {\rm sgn} \left( \frac{\det(I - A_L)
\frac{\partial}{\partial \eta} \det(I + A_L)}
{\det(I - A_L A_R)} \Big|_{(0,0)} \right) h_2(\hat{\mu})$.
\label{th:c2pd}
\end{theorem}

A proof is given in Appendix \ref{sec:PROOFND}.
Conditions (\ref{it:pdSing2})-(\ref{it:pdTrans2})
are analogous to the first three conditions of Theorem \ref{th:c21dpd}.
Condition (\ref{it:pdCenterMan2}) is a non-degeneracy condition
that ensures the period-doubled cycles collide with the \sw~in the generic manner.
If this condition is not satisfied, a higher codimension
situation that involves linear separation between $h_1(\mu)$
and $h_2(\mu)$ may arise.
Finally condition (\ref{it:pdNDRHM2}) is necessary to ensure
that non-degenerate, two-cycles are created at $\mu = 0$.

Theorem \ref{th:c2pd} essentially predicts the same
basic bifurcation structure as Theorem \ref{th:c21dpd}.
The relative position of $h_1$ and $h_2$
and the criticality of the period-doubling bifurcation
are determined by the sign of $c_0$.
In one dimension, $c_0$ is a simple function
of coefficients of nonlinear terms of the left half-map (\ref{eq:c0PD}).
To obtain the value of $c_0$ in higher dimensions,
one may derive an expression
for the restriction of the left half-map to the center manifold,
then apply the one dimensional result.
This is done for an example below.
An explicit expression for $c_0$ in $N$ dimensions
is probably too complicated to be of practical use.

To illustrate Theorem \ref{th:c2pd} consider the following map
on $\mathbb{R}^2$ using $x = (s,y)$:
\begin{equation}
\begin{split}
s' & = -\frac{1}{2} \mu - \frac{1}{2} s + y \;, \\
y' & = \frac{1}{3} s - \frac{1}{6} |s| - \frac{3}{2} \eta s + \frac{1}{4} s^2 \;.
\end{split}
\label{eq:pdMapEx}
\end{equation}
This map is easily rewritten in the piecewise form (\ref{eq:pwsMap}) with
(\ref{eq:fiForm}) and as we will show
(\ref{eq:pdMapEx}) satisfies the assumptions of Theorem \ref{th:c2pd}.
When $\mu = 0$ the origin is a fixed point of (\ref{eq:pdMapEx})
which lies on the \sw, $s = 0$.
Consequently the $\eta$-axis corresponds to a locus of border-collision bifurcations. 

\begin{figure}[ht]
\begin{center}
\includegraphics[width=10.8cm,height=9cm]{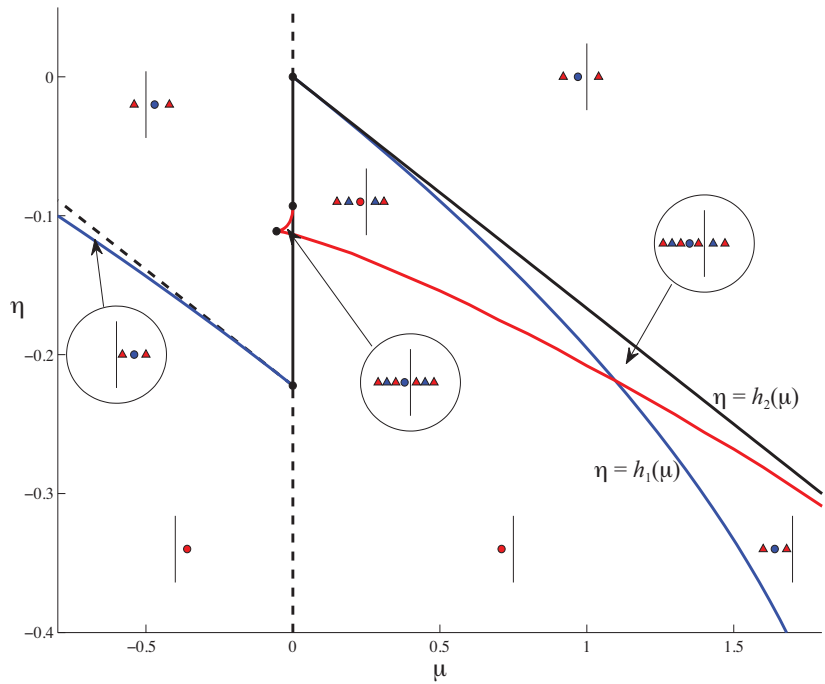}
\caption{
A bifurcation set of the two-dimensional,
piecewise-$C^\infty$ map (\ref{eq:pdMapEx}).
Border-collision bifurcations of fixed points and two-cycles are indicated by black
curves which are dashed when no bifurcation occurs in a topological sense.
Blue curves correspond to period-doubling bifurcations of a fixed point;
red curves correspond to saddle-node bifurcations of a two-cycle.
As in Fig.~\ref{fig:bsPD},
schematics showing fixed points and two-cycles are included.
Dynamics local to the codimension-two points
$(\mu,\eta) = (0,0)$ and $(0,-\frac{2}{9})$ are described by Theorem \ref{th:c2pd}.
\label{fig:bsPDex}
}
\end{center}
\end{figure}

Fig.~\ref{fig:bsPDex} shows a two-dimensional bifurcation diagram
of (\ref{eq:pdMapEx}) near $\mu = \eta = 0$.
Period-doubling bifurcations occur along $\eta = h_1(\mu)$
and the period-doubled solutions undergo border-collision along $\eta = h_2(\mu)$;
these curves are given explicitly by
\begin{equation}
\begin{split}
h_1(\mu) & = -\frac{1}{6} \mu - \frac{1}{48} \mu^2 + O(\mu^3) \;, \\
h_2(\mu) & = -\frac{1}{6} \mu \;,
\end{split}
\label{eq:h1h2ex}
\end{equation}
(actually $h_1^{-1}(\eta) \equiv -6 \eta - \frac{9}{2} \eta^2$).
The curves are tangent at $\mu = \eta = 0$ as predicted by Theorem \ref{th:c2pd}.
Near this point the dynamics resemble Fig.~\ref{fig:bsPD}C
although there is an additional unstable two-cycle present
that does not undergo bifurcation at $\mu = \eta = 0$.
This two-cycle collides with the stable two-cycle born at $\mu = 0$ for small $\eta < 0$
along a locus of saddle-node bifurcations (shown red in Fig.~\ref{fig:bsPDex}).
The saddle-node locus emanates from the $\eta$-axis
at $\eta = \frac{-4+\sqrt{10}}{9} \approx -0.093$ which is where
the stable two-cycle loses stability.
Also a cusp bifurcation of period-two orbits
occurs at $(\mu,\eta) = (-\frac{1}{18},-\frac{1}{9})$.

A second simultaneous border-collision and period-doubling 
bifurcation occurs on the $\eta$-axis at $\eta = -\frac{2}{9}$.
As predicted by Theorem \ref{th:c2pd} a locus of period-doubling bifurcations
and border-collision bifurcations of the period-doubled solution
emanate from the codimension-two point and are tangent to one another here.
Locally the bifurcation set resembles Fig.~\ref{fig:bsPD}A.

We now perform calculations to verify the details of
Theorem \ref{th:c2pd} for the map (\ref{eq:pdMapEx}) near $\mu = \eta = 0$.
Here $\varrho^{\sf T}(0,0) = [1,1]$, $b(0,0) = [-\frac{1}{2},0]^{\sf T}$
and $\det(I - A_L(0,0)) = 1$, thus $\hat{\mu} = -\mu$
and therefore $h_1(\mu)$ corresponds to admissible solutions when $\mu > 0$,
matching Fig.~\ref{fig:bsPDex}.
Also, $\frac{\partial}{\partial \eta} \det(I + A_L) \Big|_{(0,0)} = \frac{3}{2}$
and $\det(I - A_L(0,0) A_R(0,0)) = \frac{1}{6}$,
hence by the last statement of Theorem \ref{th:c2pd}, two-cycles
born on $h_2$ with points
on each side of the \sw~are admissible below $h_2$,
in agreement with Fig.~\ref{fig:bsPDex}.

Using a series expansion (see also the proof of Theorem \ref{th:c2pd}),
one finds that the left half-map of (\ref{eq:pdMapEx})
restricted to the center manifold is given by
\begin{equation}
s' = \hat{\mu} - s + \frac{1}{2} s^2 - \frac{2}{3} \hat{\mu} s + \eta s
+ \frac{1}{3} \hat{\mu}^2 - 2 \hat{\mu} \eta - \frac{1}{3} s^3 + \cdots
\nonumber
\end{equation}
Here the coefficient $-( \frac{\partial a_L}{\partial \hat{\mu}} + p) \Big|_{(0,0)}$,
given in Theorem \ref{th:c21dpd}, is equal to $\frac{1}{6}$.
Also $c_0 = -\frac{1}{12}$.
Thus $h_1(\mu) = -\frac{1}{6} \mu + O(\mu^2)$ and
$h_2(\mu) - h_1(\mu) = \frac{1}{48} \mu^2 + O(\mu^3)$
matching the explicit results given above.

\section{Conclusions}
\label{sec:CONC}

A two-parameter bifurcation diagram
provides a useful and concise picture of the behavior of a dynamical system
near a codimension-two point.
Loci or curves in such diagrams correspond to codimension-one phenomena.
The intersection of two or more bifurcation curves at a codimension-two point
occurs in a manner that is generally determined by the 
corresponding codimension-one bifurcations.
We have shown that for a general \pws, continuous map,
loci of border-collision and period-doubling bifurcations
of a single fixed point generically intersect non-tangentially
(unlike the analogous intersection of a \bcb~and a saddle-node bifurcation \cite{Si08}).
Moreover the period-doubled solution undergoes border-collision along a locus
that emanates from the codimension-two point and is tangent to the 
period-doubling curve here.

We have determined all local dynamical phenomena that occurs
when this scenario arises generically in a one-dimensional system.
In particular we have shown exactly when chaotic dynamics may be generated.
In higher dimensions \bcb s near the codimension-two point may
in addition generate quasiperiodic solutions such as invariant circles
or other more complicated orbits.

\appendix

\section{Proof of Theorem \ref{th:c21dpd}}
\label{sec:PROOF1D}

The theorem is proved in five steps.
First, the function $h_1$ is calculated by finding where
$x^{*(L)}$ has an associated multiplier of $-1$.
Second, $h_2$ is calculated from an expression for $f^{(L)^{\scriptstyle 2}}$.
Third, an explicit computation of the \pws, second iterate map
allows the \bcb~of the two-cycle that occurs along $\eta = h_2(\mu)$
to be classified.
Finally, in steps four and five we prove parts (\ref{it:resncycles})
and (\ref{it:reschaos}) respectively.

Since $f^{(L)}$ (\ref{eq:1dcoeffs}) is $C^k$ and $k \ge 3$, we can write
\begin{equation}
\begin{split}
a_L(\mu,\eta) & = -1 + \alpha_1 \mu + \eta + \alpha_3 \mu^2
+ \alpha_4 \mu \eta + \alpha_5 \eta^2 + o(2) \;, \\
b(\mu,\eta) & = 1 + \beta_1 \mu + \beta_2 \eta + O(2) \;, \\
p(\mu,\eta) & = \gamma_0 + \gamma_1 \mu + \gamma_2 \eta + o(1) \;, \\
q(\mu,\eta) & = \delta_0 + o(0) \;,
\end{split}
\label{eq:coeffsPD}
\end{equation}
where we have used hypotheses of the theorem to simplify the coefficients.

{\bf Step 1:}~~{\em Compute the function $h_1$ which corresponds to
the existence of a fixed point of the left half-map with multiplier $-1$.}\\
By the implicit function theorem the left half-map has the fixed point
\begin{equation}
x^{*(L)}(\mu,\eta) = \frac{1}{2} \mu + k_1 \mu^2 + k_2 \mu \eta + k_3 \eta^2 + O(3) \;,
\label{eq:xStarPD}
\end{equation}
which is a $C^k$ function of $\mu$ and $\eta$ and locally satisfies
$f^{(L)}(x^{*(L)}(\mu,\eta);\mu,\eta) - x^{*(L)}(\mu,\eta) \equiv 0$.
By a second order expansion, it is determined that
the scalar coefficients, $k_i$, have the values
\begin{equation}
\begin{split}
k_1 & = \frac{\beta_1}{2} + \frac{\alpha_1}{4} +
\frac{\gamma_0}{8} \;, \\
k_2 & = \frac{\beta_2}{2} + \frac{1}{4} \;, \\
k_3 & = 0 \;.
\end{split}
\label{eq:kiPD}
\end{equation}
The multiplier associated with $x^{*(L)}(\mu,\eta)$ is found by
substituting (\ref{eq:xStarPD}) into $D_x f^{(L)}(x;\mu,\eta)$ (see (\ref{eq:1dcoeffs})).
We obtain
\begin{eqnarray}
D_x f^{(L)}(x^{*(L)}(\mu,\eta);\mu,\eta) & = & -1 + (\alpha_1 + \gamma_0) \mu + \eta
+ \left( \alpha_3 + 2 k_1 \gamma_0 + \gamma_1 +
\frac{3 \delta_0}{4} \right) \mu^2
\nonumber \\
& & +~(\alpha_4 + 2 k_2 \gamma_0 + \gamma_2) \mu \eta
+ \alpha_5 \eta^2 + o(2) \;.
\label{eq:xStarMultPD}
\end{eqnarray}
The implicit function theorem implies that
\begin{equation}
h_1(\mu) = -(\alpha_1 + \gamma_0) \mu + l_1 \mu^2 + o(\mu^2) \;,
\label{eq:h1proofPD}
\end{equation}
is a $C^{k-1}$ function that locally satisfies
$D_x f^{(L)}(x^{*(L)}(\mu,h_1(\mu));\mu,h_1(\mu)) + 1 \equiv 0$,
where $l_1$ is found by substituting
(\ref{eq:h1proofPD}) into (\ref{eq:xStarMultPD}),
\begin{equation}
l_1 = -\left( \alpha_3 + 2 k_1 \gamma_0 + \gamma_1
+ \frac{3 \delta_0}{4} \right)
+ (\alpha_4 + 2 k_2 \gamma_0 + \gamma_2)(\alpha_1 + \gamma_0)
- \alpha_5 (\alpha_1 + \gamma_0)^2 \;.
\label{eq:l1PD}
\end{equation}
When $c_0 \ne 0$, the existence of period-doubling bifurcations along $\eta = h_1(\mu)$
for small $\mu$ is verified by checking the three conditions of the
standard theorem, see for instance \cite{GuHo86}:
\begin{enumerate}[label=\roman{*}),ref=\roman{*}]
\item
(singularity)
by construction, $D_x f^{(L)}(x^{*(L)}(\mu,h_1(\mu));\mu,h_1(\mu)) \equiv -1$,
\item
(transversality)
$\left( \frac{\partial f^{(L)}}{\partial \eta}
\frac{\partial^2 f^{(L)}}{\partial x^2}
+ 2 \frac{\partial^2 f^{(L)}}{\partial x \partial \eta} \right)
\Big|_{\eta = h_1(\mu)} = 2 + O(\mu) \ne 0$,
\item
(non-degeneracy)
$\left( \frac{1}{2} \left( \frac{\partial^2 f^{(L)}}{\partial x^2} \right)^2
+ \frac{1}{3} \frac{\partial^3 f^{(L)}}{\partial x^3} \right)
\bigg|_{\eta = h_1(\mu)} = 2 c_0 + O(\mu) \ne 0$.
\end{enumerate}
Consequently we have proved (\ref{it:resh1}) of the theorem.

{\bf Step 2:}~~{\em Compute the function $h_2$ which corresponds to the existence of
a period-two orbit of the left half-map at $x = 0$.}\\
To determine $h_2$, we compute the $C^k$ function
$f^{(L)^{\scriptstyle 2}}(0;\mu,\eta)
= f^{(L)}(\mu b(\mu,\eta);\mu,\eta)$ which may be written as
$f^{(L)^{\scriptstyle 2}}(0;\mu,\eta) = \mu g(\mu,\eta)$ where
\begin{equation}
g(\mu,\eta) = b(\mu,\eta)
\Big( 1 + a_L(\mu,\eta) + \mu b(\mu,\eta) p(\mu,\eta)
+ \mu^2 b^2(\mu,\eta) q(\mu,\eta) \Big) + o(2) \;,
\label{eq:gPD}
\end{equation}
is $C^{k-1}$.
Substitution of (\ref{eq:coeffsPD}) into (\ref{eq:gPD}) produces
\begin{eqnarray}
g(\mu,\eta) & = & (\alpha_1 + \gamma_0) \mu + \eta
+ (\alpha_1 \beta_1 + \alpha_3 + 2 \beta_1 \gamma_0
+ \gamma_1 + \delta_0) \mu^2 \nonumber \\
& & +~(\beta_1 + \alpha_1 \beta_2 + \alpha_4 + 2 \beta_2 \gamma_0
+ \gamma_2) \mu \eta + (\beta_2 + \alpha_5) \eta^2 + o(2) \;.
\hspace{15mm}
\label{eq:g2PD}
\end{eqnarray}
The implicit function theorem implies that
\begin{equation}
h_2(\mu) = -(\alpha_1 + \gamma_0) \mu
+ l_2 \mu^2 + o(\mu^2) \;,
\label{eq:h2PD}
\end{equation}
is a $C^{k-1}$ function that locally satisfies $g(\mu,h_2(\mu)) \equiv 0$.
By substituting (\ref{eq:h2PD}) into (\ref{eq:g2PD}) we deduce
\begin{eqnarray}
l_2 & = & -(\alpha_1 \beta_1 + \alpha_3 + 2 \beta_1 \gamma_0 + \gamma_1 + \delta_0)
+ (\beta_1 + \alpha_1 \beta_2 + \alpha_4 + 2 \beta_2 \gamma_0 + \gamma_2)
(\alpha_1 + \gamma_0) \nonumber \\
& & -~(\beta_2 + \alpha_5)
(\alpha_1 + \gamma_0)^2 \;.
\label{eq:l2PD}
\end{eqnarray}
Subtracting (\ref{eq:l1PD}) from (\ref{eq:l2PD}) produces
(after algebraic simplification)
\begin{equation}
l_2 - l_1 = -\frac{\gamma_0^2 + \delta_0}{4} \;, \nonumber
\end{equation}
which proves (\ref{eq:h2PDth}) in the statement of the theorem.

{\bf Step 3:}~~{\em Determine all two-cycles to verify the phase portraits
in Fig.~\ref{fig:bsPD}.}\\
For small, fixed $\mu < 0$, the two-cycle generated in a period-doubling bifurcation at
$\eta = h_1(\mu)$ undergoes a \bcb~when it collides with the \sw~at $\eta = h_2(\mu)$.
The stability and relative admissibility of two-cycles
emanating from this \bcb~is found by determining
a map of the form (\ref{eq:1dmap}) that describes the bifurcation.
Such a map may be obtained by computing the second iterate of (\ref{eq:1dmapPD})
and replacing $\eta$ with the new parameter
\begin{equation}
\hat{\eta} = \eta - h_2(\mu) \;, \nonumber
\end{equation}
which controls the border-collision.
Near $\hat{\eta} = 0$, when $\mu < 0$
period-two orbits are guaranteed to have one negative point, so we compute
$f^{(L)} \circ f^{(L)}$ and $f^{(L)} \circ f^{(R)}$ which leads to
\begin{equation}
f^2(x;\mu,\hat{\eta}) = \left\{ \begin{array}{lc}
(\mu + O(\mu^2)) \hat{\eta} +
(1 - c_0 \mu^2 + o(\mu^2)) x + O(|x,\hat{\eta}|^2), & x \le 0 \\
(\mu + O(\mu^2)) \hat{\eta} +
(-a_0^{(R)} + O(\mu)) x + O(|x,\hat{\eta}|^2), & x \ge 0
\end{array} \right.
\label{eq:f2proofPD}
\end{equation}
A fixed point, $x^{*(LL)}$, of the left half-map of (\ref{eq:f2proofPD})
corresponds to a two-cycle of (\ref{eq:1dmapPD}) with both points
negative-valued, when admissible.
By (\ref{eq:f2proofPD}), if $c_0 < 0$ [$c_0 > 0$], then
$x^{*(LL)}$ is unstable and admissible when $\hat{\eta} < 0$
[stable and admissible when $\hat{\eta} > 0$].
Similarly, a fixed point, $x^{*(RL)}$ of the right half-map of (\ref{eq:f2proofPD})
corresponds to a two-cycle of (\ref{eq:1dmapPD}) comprised of
two points with different signs, when admissible.
If $|a_0^{(R)}| < 1$, then $x^{*(RL)}$ is stable, otherwise it is unstable.
Furthermore if $a_0^{(R)} < -1$ [$a_0^{(R)} > -1$],
then $x^{*(RL)}$ is admissible when $\hat{\eta} > 0$ [$\hat{\eta} < 0$].
Lastly, since $a_0^{(R)} \ne -1$, the fixed point, $x^{*(RL)}$,
is unique for small $\mu$ and $\eta$.
Therefore the two-cycle created at $\mu = 0$
collides with the \sw~at $\eta = h_2(\mu)$.
These statements verify all two-cycles shown in Fig.~\ref{fig:bsPD}.

{\bf Step 4:}~~{\em Show that $f$, (\ref{eq:1dmapPD}),
has no $n$-cycles with $n \ge 3$, verifying part (\ref{it:resncycles}).}\\
Except when $a_0^{(R)} = 0$ our proof is founded on the knowledge
that a one-dimensional map cannot have an $n$-cycle with $n \ge 3$
contained in an interval
on which the map is monotone.
For the special case $a_0^{(R)} = 0$ we used several additional logical arguments
including the contraction mapping theorem.



{\bf Case I:}~~$a_0^{(R)} < 0$.\\
In this case there exists $\delta > 0$ such that
$\forall \mu, \eta \in [-\delta,\delta]$,
$f$ is decreasing on $[-\delta,\delta]$.
Consequently $f$ has no $n$-cycles with $n \ge 3$ on $[-\delta,\delta]$.

{\bf Case II:}~~$a_0^{(R)} = 0$.\\
First we construct an interval containing the origin on which $f$
is forward invariant.
Since the two components of $f$ are differentiable,
there exists $\hat{\delta} > 0$ such that
$\forall \mu, \eta, x \in [-\hat{\delta},\hat{\delta}]$,
\begin{equation}
\left| \frac{\partial f^{(R)}}{\partial x} \right|,
\left| \frac{\partial f^{(L)}}{\partial x} + 1 \right|,
\left| b(\mu,\eta) - 1 \right| \le \frac{1}{3} \;.
\nonumber
\end{equation}
Let $\delta = \frac{1}{8} \hat{\delta}$,
let $I = [-4 \delta, 8 \delta]$
and assume $\mu,\eta \in [-\delta,\delta]$.
Then $f(I) \subset I$ because for any $x \in I$,
we may assume $x > 0$ when calculating a lower bound for $f(x)$:
$f(x) \ge \mu b - \frac{1}{3} x \ge -\frac{4}{3} \delta - \frac{8}{3} \delta = -4 \delta$,
and we may assume $x < 0$ when calculating an upper bound:
$f(x) \le \mu b - \frac{4}{3} x \le \frac{4}{3} \delta + \frac{4}{3} 4 \delta < 8 \delta$.

If $\mu \ge 0$ then the image of any $x < 0$ under $f$ is positive
and since the right half-map is contracting and the slope of the
left half-map is near $-1$,
$f^2$ is a contraction on $I$,
indeed $|D_x f^2(x)| \le \frac{4}{3} \frac{1}{3} < 1$.
By the contraction mapping theorem, iterates of $f^2$ approach a unique
fixed point, which must be $x^{*(R)}$, therefore
$f$ cannot have an $n$-cycle with $n \ge 3$ on $\mathcal{N}$ (\ref{eq:nbd}).

\begin{figure}[ht]
\begin{center}
\includegraphics[width=12cm,height=5cm]{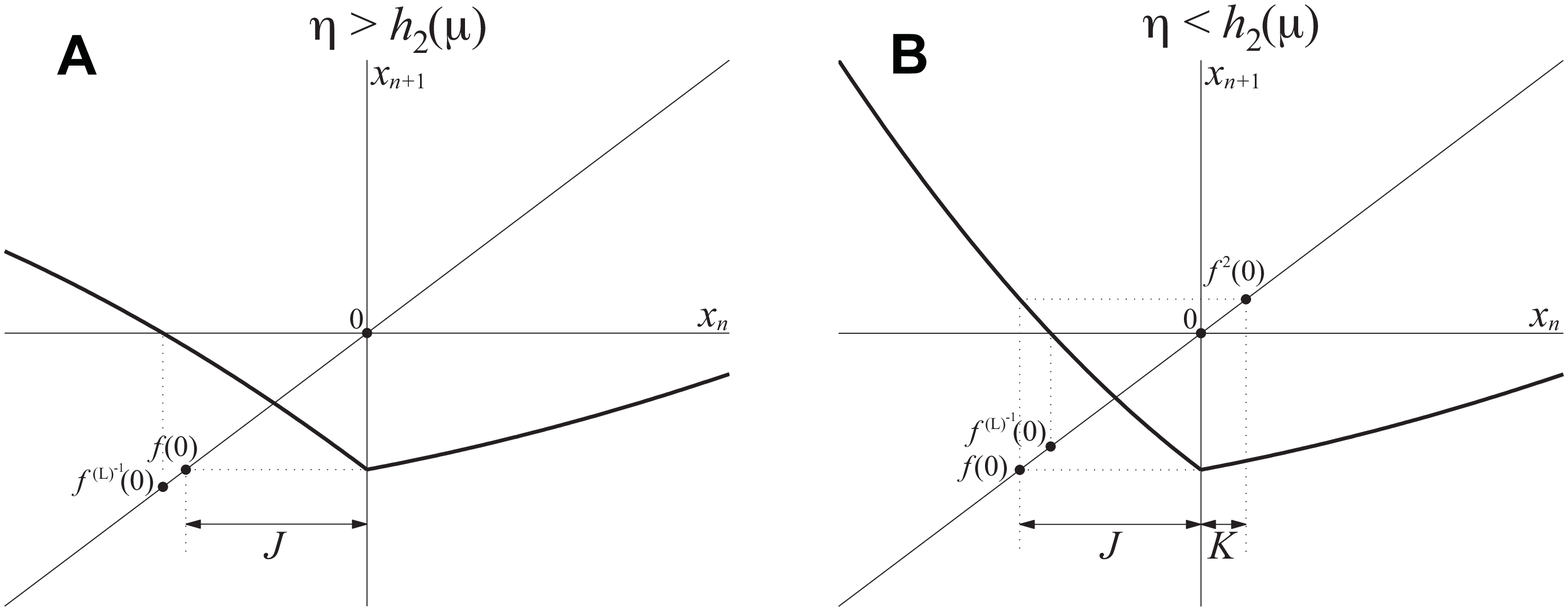}
\caption{
Schematic diagrams illustrating the map $f$, (\ref{eq:1dmapPD}),
when $\mu < 0$ and $\eta > h_2(\mu)$
in panel A and $\eta < h_2(\mu)$ in panel B.
When $\eta = h_2(\mu)$, $f^{(L)^{\scriptstyle -1}}(0) = f(0)$.
\label{fig:a0Rgen}
}
\end{center}
\end{figure}

So we now assume $\mu < 0$.
Recall that when $\eta = h_2(\mu)$, $f^{(L)^{\scriptstyle 2}}(0) = 0$,
i.e.~$f(0) = f^{(L)^{\scriptstyle -1}}(0)$,
and when $\eta > h_2(\mu)$, $f(0) > f^{(L)^{\scriptstyle -1}}(0)$.
Therefore if $\eta \ge h_2(\mu)$, see Fig.~\ref{fig:a0Rgen}-A,
the interval $[f^{(L)^{\scriptstyle -1}}(0),0]$ is forward invariant.
Since $f$ is monotone on the interval
\begin{equation}
J = [f(0),0] \;,
\label{eq:JPD}
\end{equation}
it contains no $n$-cycles with $n \ge 3$ here.
On the other hand if $\eta < h_2(\mu)$ as in Fig.~\ref{fig:a0Rgen}-B,
the forward orbit of any point in $[f^{(L)^{\scriptstyle -1}}(0),0]$
is either contained in this interval (where $f$ is monotone) or
enters $[f(0),f^{(L)^{\scriptstyle -1}}(0)]$
and approaches the stable two-cycle
(this may be shown formally by the contraction mapping theorem).
In any case, an $n$-cycle of $f$ with $n \ge 3$ cannot
include a point in $J = [f(0),0]$.

Suppose for a contradiction $f$ has an $n$-cycle with $n \ge 3$,
call it $\{ x_0 \ldots x_{n-1} \} \subset [-\delta,\delta] \setminus J$.
We now show that the restriction of $f^2$ to this $n$-cycle is a contraction
giving a contradiction according to the
contraction mapping theorem.
Specifically we will argue that for any $i$ and $j$,
\begin{equation}
|f^2(x_j) - f^2(x_i)| < \frac{4}{9} |x_j - x_i| \;.
\label{eq:contract49}
\end{equation}
There are several cases to consider.
For example, suppose that $x_j > 0$, $x_i < 0$ and $f^2(x_j) < f^2(x_i)$.
Then $f(x_j) > f(0) - \frac{1}{3} x_j$ thus
$f^2(x_j) > f(0) + \frac{1}{3} f(0) - \frac{1}{9} x_j$.
Also $f(x_i) < f(0) - \frac{4}{3} x_i$ thus
$f^2(x_i) < f(0) + \frac{1}{3} f(0) - \frac{4}{9} x_i$.
Subtracting these inequalities leads to (\ref{eq:contract49}).
The remaining cases can be shown similarly.

{\bf Case III:}~~$a_0^{(R)} > 0$.\\
{\em Case IIIa:}
First suppose that $\mu \ge 0$.
There exists $\delta_1 > 0$ such that whenever $\mu \in [0,\delta_1]$
and $\eta \in [-\delta_1,\delta_1]$,
$f$ is decreasing on $[-\delta_1,0]$ and increasing on $[0,\delta_1]$.
Since $\mu \ge 0$, $\forall x \in [-\delta_1,\delta_1]$, $f(x) \ge 0$.
Thus any periodic solution of $f$ in $[-\delta_1,\delta_1]$
must lie entirely in $[0,\delta_1]$.
But $f$ is increasing on this interval thus has no $n$-cycles with $n \ge 3$ here.\\
{\em Case IIIb:}
Now suppose $\mu < 0$ and $\eta \ge h_2(\mu)$.
There exists $\delta_2 > 0$ such that whenever $\mu \in [-\delta_2,0)$
and $\eta \in [h_2(\mu),\delta_2]$,
$f$ is decreasing on $[-\delta_2,0]$ and increasing on $[0,\delta_2]$,
see Fig.~\ref{fig:a0Rgen}-A.
Suppose for a contradiction, $f$ has an $n$-cycle with $n \ge 3$ on $[-\delta_2,\delta_2]$.
Such an $n$-cycle must enter both $[-\delta_2,0)$ and $(0,\delta_2]$, and so
it includes a point $x \in (0,\delta_2]$
with $f(x) \in [-\delta_2,0)$.
But $f(x)$ lies in $J$, (\ref{eq:JPD}), and
since $\eta \ge h_2(\mu)$, $f^2(0) \le 0$, hence $f(J) \subset J$.
Thus the forward orbit of $x$ cannot return to $(0,\delta_2]$
which contradicts the assumption that $x$ belongs to an $n$-cycle.\\
{\em Case IIIc:}
Finally suppose $\mu < 0$, $\eta < h_2(\mu)$ and $0 < a_0^{(R)} < 1$.
Then there exists $\delta_3 > 0$ such that
whenever $\mu \in [-\delta_3,0)$,
$\eta \in [-\delta_3,h_2(\mu))$
and $x \in [-\delta_3,\delta_3]$,
\begin{eqnarray*}
-\left(1 + \frac{1-a_0^{(R)}}{2(1+a_0^{(R)})}\right)
\le & f^{(L)'}(x) & \le -\frac{1}{2} \;, \\
\frac{a_0^{(R)}}{2} \le & f^{(R)'}(x) & \le \frac{1+a_0^{(R)}}{2} \;,
\end{eqnarray*}
where the particular non-symmetric $a_0^{(R)}$-dependent bounds on the slopes
have been chosen to simplify the subsequent analysis.
Suppose for a contradiction $f$ has an $n$-cycle with $n \ge 3$ in $[-\delta_3,\delta_3]$.
As before, such an $n$-cycle must include a point $x \in (0,\delta_3]$
with $f(x) \in [-\delta_3,0)$.
Again $f(x) \in J$, but here $\eta < h_2(\mu)$,
thus the forward orbit of $x$ may return to the right of the origin
but must enter the interval $K = [0,f^2(0)]$, see Fig.~\ref{fig:a0Rgen}-B.
Notice $f^3(0) \le f(0) + \frac{1+a_0^{(R)}}{2} f^2(0)$.
Then by evaluating the first-order expansion of $f^{(L)}$ about $f(0)$
at $f^3(0)$,
we arrive at
$f^4(0) \ge f^2(0) - \left( 1 + \frac{1}{2} \frac{1-a_0^{(R)}}{1+a_0^{(R)}} \right)
\frac{1+a_0^{(R)}}{2} f^2(0) > 0$, therefore $f^2(K) \subset K$.
Then $\forall y \in K$,
$|f^{2{\textstyle \hspace{.3mm}'}}(y)| \le
\left( 1 + \frac{1}{2} \frac{1-a_0^{(R)}}{1+a_0^{(R)}} \right) 
\frac{1+a_0^{(R)}}{2} < 1$.
Thus $f^2 : K \to K$ is a contraction mapping and therefore iterates of $f$
that enter $K$ cannot belong to an $n$-cycle with $n \ge 3$.\\
Finally let $\delta = {\rm min}(\delta_1,\delta_2,\delta_3)$,
then the result, (\ref{eq:1dmapPD}), is proved.

\begin{figure}[ht]
\begin{center}
\includegraphics[width=7.2cm,height=6cm]{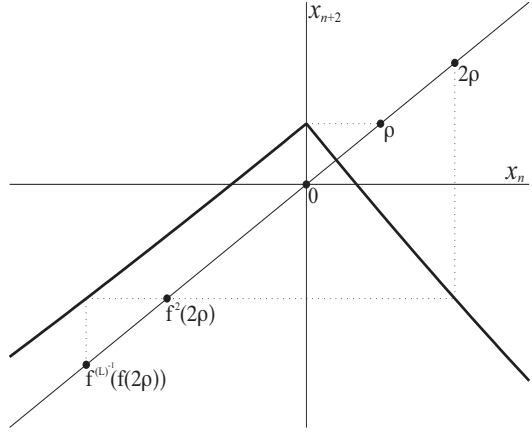}
\caption{
The second iterate of (\ref{eq:pwsMap}), $f^2$,
near $x = 0$ when $\mu < 0$, $\eta < h_2(\mu)$ and $a_0^{(R)} > 1$.
The slope of this map is approximately unity left of the origin
and approximately $-a_0^{(R)}$ right of the origin.
As a consequence, see the text, the map exhibits chaos
in $T = \left[ f^{(L)^{\scriptstyle -1}}(f(2\rho)), 2\rho \right]$.
\label{fig:f2PDschem}
}
\end{center}
\end{figure}

{\bf Step 5:}~~{\em Prove the existence of chaos when
$\mu < 0$, $\eta < h_2(\mu)$ and $a_0^{(R)} > 1$.}\\
We first construct a trapping set, $T$, for $f^2$.
Let $\rho(\mu,\eta) = f^2(0;\mu,\eta)$.
Since $\mu < 0$, by (\ref{eq:f2proofPD}), $\rho(\mu,\eta) = O(2)$.
Since $\eta < h_2(\mu)$, $\rho > 0$.
Therefore
\begin{equation}
f(2\rho) = f(0) + (2 a_0^{(R)} + O(1)) \rho + O(\rho^2) < 0 \;.
\nonumber
\end{equation}
Consequently,
\begin{eqnarray}
f^2(2\rho) & = & (1 - 2 a_0^{(R)} + O(1)) \rho + O(\rho^2)
\label{eq:f22rho}
\;, \\
f^{(L)^{\scriptstyle -1}}(f(2\rho)) & = & (-2 a_0^{(R)} + O(1)) \rho + O(\rho^2)
\label{eq:Tmin}
\;.
\end{eqnarray}
Let
\begin{equation}
T = \left[ f^{(L)^{\scriptstyle -1}}(f(2\rho)), 2\rho \right] \;,
\nonumber
\end{equation}
It follows, see Fig.~\ref{fig:f2PDschem},
that $\forall x \in T$,
$f^2(2\rho) \le f^2(x) \le \rho$.
But $\rho, f^2(2\rho) \in {\rm int}(T)$,
thus $f^2(T) \subset {\rm int}(T)$,
i.e.~$T$ is a trapping set for $f^2$.
Furthermore,
$\bigcap_{i=0}^\infty f^{2i}(T)$,
is an attracting set for $f^2$.

Let $M \in \mathbb{Z}$ with $M \ge 2 a_0^{(R)} + 2$.
We now show that $\forall x \in T$, $\exists j \in \mathbb{Z}$
with $0 \le j < M$ such that $f^{2j}(x) > 0$.
Suppose otherwise.
Then whenever $0 \le j < M$,
$f^{2j}(x) \le 0$ and thus
\begin{equation}
f^{2j}(x) = (j + O(1)) \rho + (1 + O(1)) x + O(x^2) \;.
\label{eq:f2jx}
\end{equation}
Since $f^2$ is increasing on $T_-$,
it suffices to consider $x = f^{(L)^{\scriptstyle -1}}(f(2\rho))$.
By combining (\ref{eq:Tmin}) and (\ref{eq:f2jx}) we find
\begin{eqnarray}
f^{2(M-1)}(x) & \ge & (2 a_0^{(R)} + 1 + O(1)) \rho
+ (1 + O(1))(-2 a_0^{(R)} + O(1)) \rho + O(\rho^2) \nonumber \\
& = & (1 + O(1)) \rho + O(\rho^2) > 0 \nonumber \;,
\end{eqnarray}
which is a contradiction.

Now consider the map $f^{2M} : T \to T$.
Note that $f^2$ has one critical point on $T$, namely 0.
Thus $y \in T$ is a critical point of $f^{2M}$ if
$f^{2j}(y) = 0$ for some $0 \le j < M$.
Consequently $f^{2M}$ has at most $2^M - 1$ critical points, $y_i$,
and $f^{2M}$ is differentiable on $T$ except at each $y_i$.
For all $x \in T$, $x \ne y_i$,
\begin{equation}
f^{2M{\textstyle \hspace{.3mm}'}}(x) =
\prod_{j=0}^{M-1} f^{2{\textstyle \hspace{.3mm}'}}(f^{2j}(x)) \nonumber \;.
\end{equation}
Let $m \ge 1$ be the number of iterates, $f^{2j}(x)$, that are positive.
Since left of zero the slope of $f^2$ is near 1
and right of zero the slope is near $-a_0^{(R)}$,
we have
\begin{equation}
\left| f^{2M{\textstyle \hspace{.3mm}'}}(x) \right| =
a_0^{(R)^{\scriptstyle m}} + O(1) > \frac{1 + a_0^{(R)}}{2} > 1 \;,
\nonumber
\end{equation}
for any $x \ne y_i$, for sufficiently small $\mu, \eta$.
Consequently $f^{2M} : T \to T$ is {\em piecewise expanding} \cite{Ro04}.
By the theorem of Li and Yorke \cite{LiYo78},
$f^{2M}$ is chaotic on $T$.
Thus $f$ is chaotic on $T \cup f(T)$.
\qed

\section{Proof of Theorem \ref{th:c2pd}}
\label{sec:PROOFND}

The theorem is proved in three steps.
In the first two steps a center manifold, $W^c$, for the left half-map, $f^{(L)}$,
is constructed and the lowest order terms of the restriction of $f^{(L)}$
to $W^c$ are calculated.
In the third step the two-cycle that has points in each half-plane
is computed explicitly.

{\bf Step 1:}~~{\em Compute $W^c$.}\\
It is convenient to extend the map to $\mathbb{R}^N \times \mathbb{R}^2$
by including the parameters in the extended phase space,
\begin{equation}
F = \left[ \begin{array}{c}
x' \\ \mu' \\ \eta'
\end{array} \right]
= \left[ \begin{array}{c}
f^{(L)}(x;\mu,\eta) \\ \mu \\ \eta
\end{array} \right] \;.
\nonumber
\end{equation}
Then
\begin{equation}
D F(0;0,0) = \left[ \begin{array}{ccccc}
A_L(0,0) & \Bigg| & b(0,0) & \Bigg| & 0 \\
\hline
\begin{array}{c} 0 \\ 0 \end{array} & \Bigg| &
\begin{array}{c} 1 \\ 0 \end{array} & \Bigg| &
\begin{array}{c} 0 \\ 1 \end{array}
\end{array} \right] \;,
\nonumber
\end{equation}
has a three-dimensional center space
\begin{equation}
E^c = {\rm span} \left\{
\left[ \begin{array}{c} v(0,0) \\ 0 \\ 0 \end{array} \right],
\left[ \begin{array}{c} \varphi \\ 1 \\ 0 \end{array} \right],
e_{N+2} \right\} \;,
\nonumber
\end{equation}
where
\begin{equation}
\varphi = (I - A_L(0,0))^{-1} b(0,0) \;.
\label{eq:varphiProofPD}
\end{equation}
Notice
\begin{equation}
e_1^{\sf T} \varphi =
\frac{\varrho^{\sf T}(0,0) b(0,0)}{\det(I - A_L(0,0))} \ne 0 \;.
\label{eq:varphiProofPD1}
\end{equation}
Since $e_1^{\sf T} v(0,0) = 1 \ne 0$, the center manifold, $W^c$,
may be expressed locally in terms of $s$, $\mu$ and $\eta$.
By the center manifold theorem,
there exists a $C^{k-1}$ function, $H$,
such that on $W^c$
\begin{equation}
x = H(s;\mu,\eta) = s v(0,0) + \mu \zeta + O(2)
\label{eq:HPD}
\end{equation}
where
\begin{equation}
\zeta = \varphi - (e_1^{\sf T} \varphi) v(0,0) \;.
\label{eq:zetaProofPD}
\end{equation}

{\bf Step 2:}~~{\em Determine an expression for $f^{(L)}$ on $W^c$
and verify conditions (\ref{it:pdSing})-(\ref{it:pdTrans}) of Theorem \ref{th:c21dpd}.}\\
On $W^c$
\begin{eqnarray}
s' & = & e_1^{\sf T} f^{(L)}(H(s;\mu,\eta);\mu,\eta) \nonumber \\
& = & e_1^{\sf T} b(\mu,\eta) \mu
+ e_1^{\sf T} A_L(\mu,\eta) H(s;\mu,\eta)
+ e_1^{\sf T} g^{(L)}(H(s;\mu,\eta);\mu,\eta) \;, \hspace{15mm}
\label{eq:spPD}
\end{eqnarray}
where $g^{(L)}$ denotes all terms of $f^{(L)}$ that are nonlinear in $x$.
With the given ``hatted'' variables (\ref{eq:spPD})
satisfies condition (\ref{it:pdNDSW}) of Theorem \ref{th:c21dpd} because
\begin{eqnarray}
\frac{\partial s'}{\partial \mu} \Big|_{(0;0,0)}
& = & e_1^{\sf T} b(0,0) + e_1^{\sf T} A_L(0,0) \zeta \nonumber \\
& = & e_1^{\sf T} b(0,0) + e_1^{\sf T} A_L(0,0) \varphi
- (e_1^{\sf T} \varphi) e_1^{\sf T} A_L(0,0) v(0,0)
{\rm ~~~(by~(\ref{eq:zetaProofPD}))} \nonumber \\
& = & e_1^{\sf T} \varphi + (e_1^{\sf T} \varphi) e_1^{\sf T} v(0,0)
{\rm ~~~(by~(\ref{eq:varphiProofPD}))} \nonumber \\
& = & 2 e_1^{\sf T} \varphi \;,
\nonumber
\end{eqnarray}
which with $\hat{\mu} = 2 e_1^{\sf T} \varphi \mu$,
in view of (\ref{eq:varphiProofPD1}), produces
\begin{equation}
\frac{\partial s'}{\partial \hat{\mu}} \Big|_{(0;0,0)} = 1 \;.
\nonumber
\end{equation}
A similar verification of (\ref{it:pdSing}) and (\ref{it:pdTrans})
of Theorem \ref{th:c21dpd} is now given.
By scaling $v$, we may assume that for small $\eta$,
$e_1^{\sf T} v(0,\eta) \equiv 1$.
Since (\ref{eq:HPD}) lacks a linear $\eta$ term,
\begin{equation}
H(s;0,\eta) = v(0,\eta) s + O(s^2) \;,
\nonumber
\end{equation}
and consequently
\begin{eqnarray}
s' |_{\mu = 0} & = & e_1^{\sf T} A_L(0,\eta) v(0,\eta) s + O(s^2) \nonumber \\
& = & \lambda(0,\eta) s + O(s^2) \;.
\nonumber
\end{eqnarray}
Thus
\begin{eqnarray}
\frac{\partial s'}{\partial s} \Big|_{(0;0,0)} =
\lambda(0,0) & = & -1 \;, \nonumber \\
\frac{\partial^2 s'}{\partial \eta \partial s} \Big|_{(0;0,0)} =
\frac{\partial \lambda}{\partial \eta}(0,0) & = & 1 \;,
\nonumber
\end{eqnarray}
as required.

{\bf Step 3:}~~{\em Compute the two-cycle
that has points in each half-plane to verify the final statement of the theorem.}\\
We have
\begin{equation}
f^{(L)}(f^{(R)}(x;\mu,\eta);\mu,\eta)
= (I + A_L(0,\eta)) b(0,\eta) \mu
+ A_L(0,\eta) A_R(0,\eta) x + O(|\mu,x|^2) \;.
\label{eq:fLfRproof}
\end{equation}
By applying the implicit function theorem
to (\ref{eq:fLfRproof}) in view of condition (\ref{it:pdNDRHM2}) of the theorem
we see that the two-cycle with points in each half-plane exists and
is unique for small $\mu$ and $\eta$.
When admissible, the point of this cycle in the right half-plane is
given by a $C^k$ function
\begin{equation}
x^{*(RL)}(\mu,\eta) = (I - A_L A_R)^{-1}
(I + A_L) b \Big|_{\mu = 0} \mu + O(\mu^2) \;.
\nonumber
\end{equation}
Then
\begin{eqnarray}
\frac{\partial s^{*(RL)}}{\partial \mu} \Big|_{\mu = 0}
& = & e_1^{\sf T} (I - A_L A_R)^{-1} (I + A_L) b \Big|_{\mu = 0} \nonumber \\
& = & \frac{e_1^{\sf T} {\rm adj}(I - A_L A_R) (I + A_L) b}
{\det(I - A_L A_R)} \Bigg|_{\mu = 0} \;. \nonumber
\end{eqnarray}
Since $I - A_L A_R = (I + A_L)(I - A_R) - A_L + A_R$
and $A_L$ and $A_R$ are identical in their last $N-1$ columns,
${\rm adj}(I - A_L A_R)$ and ${\rm adj}((I + A_L)(I - A_R))$
share the same first row.
Thus
\begin{eqnarray}
\frac{\partial s^{*(RL)}}{\partial \mu} \Big|_{\mu = 0}
& = & \frac{e_1^{\sf T} {\rm adj}(I - A_R) {\rm adj}(I + A_L) (I + A_L) b}
{\det(I - A_L A_R)} \Bigg|_{\mu = 0} \nonumber \\
& = & \frac{\det(I + A_L)}{\det(I - A_L A_R)} \varrho^{\sf T} b \Bigg|_{\mu = 0} \nonumber \\
& = & \frac{\frac{\partial}{\partial \eta} \det(I + A_L)}{\det(I - A_L A_R)}
\varrho^{\sf T} b \Bigg|_{(0,0)} \eta + O(\eta^2) \;. \nonumber
\end{eqnarray}
By Theorem \ref{th:c21dpd}, there exists a $C^{k-2}$ curve,
$\hat{\eta} = h_2(\hat{\mu})$, along which, $s^{*(RL)} = 0$.
Furthermore, when $\hat{\mu} \le 0$, the two-cycle
is admissible along this curve.
If $\tilde{\eta} = \eta - h_2(\hat{\mu})$, then
\begin{equation}
s^{*(RL)}(\hat{\mu},\tilde{\eta}) =
\frac{\det(I - A_L) \frac{\partial}{\partial \eta} \det(I + A_L)}
{2 \det(I - A_L A_R)} \Bigg|_{(0,0)} \hat{\mu} \tilde{\eta} + O(3) \;,
\nonumber
\end{equation}
which confirms the final statement of the theorem.
\qed


\end{document}